\documentclass[leqno,a4paper,12pt]{article}

\usepackage{mathpazo}
\usepackage{amsmath}
\usepackage{amsthm}
\usepackage{ulem}
\usepackage{amssymb}
\usepackage[T1]{fontenc}

\newtheorem{thm}{Theorem }
\newtheorem{pro}{Proposition }
\newtheorem{lem}{Lemma }

\newtheorem{rem}{Remark }

\begin{document}

\title{Dynamics of symmetric holomorphic maps on projective spaces}
\author{\textsc{Kohei Ueno}}
\date{}
\maketitle

%%%%%%%%%%%%%%%%%%%%%%%%%%%%%%%%%%%%%%%%%%%%%%%%%%%%%%%%%%%%%%%%%%%%%%%%%%%%%%%%%%%%%%%%%%
\begin{abstract}
We consider complex dynamics of a \textit{critically finite} holomorphic map 
from $\mathbf{P}^{k}$ to $\mathbf{P}^{k}$,
which has symmetries associated with the symmetric group $S_{k+2}$ acting on $\mathbf{P}^{k}$,
for each $k \ge 1$.
The Fatou set of each map of this family consists of attractive basins of
superattracting points. Each map of this family satisfies Axiom A.
\end{abstract}

%%%%%%%%%%%%%%%%%%%%%%%%%%%%%%%%%%%%%%%%%%%%%%%%%%%%%%%%%%%%%%%%%%%%%%%%%%%%%%%%%%%%%%%%%%
\section{Introduction}

For a finite group $G$ acting on $\mathbf{P}^{k}$ as projective transformations,
we say that a rational map $f$ on $\mathbf{P}^{k}$ is $G$-\hspace{0pt}\textit{equivariant}
if $f$ commutes with each element of $G$.
That is, $f \circ r=r \circ f$ for any $r \in G$,
where $\circ$ denotes the composition of maps. 
Doyle and McMullen~\cite{dm} introduced the notion of
\textit{equivariant} functions on $\mathbf{P}^{1}$ to solve quintic equations.
See also \cite{ushiki} for \textit{equivariant} functions on $\mathbf{P}^{1}$.
Crass~\cite{c1} extended Doyle and McMullen's algorithm to higher dimensions
to solve sextic equations.
Crass~\cite{c} found a good family of finite groups and \textit{equivariant} maps
for which one may say something about global dynamics.
Crass~\cite{c} conjectured that the Fatou set of each map of this family 
consists of attractive basins of superattracting points.
Although I do not know whether this family has relation to solving equations or not,
our results will give affirmative answers for the conjectures in \cite{c}.

In section 2 we shall explain an action of the symmetric group $S_{k+2}$ on $\mathbf{P}^{k}$
and properties of our $S_{k+2}$-equivariant map.
In section 3 and 4 we shall show our results
about the Fatou sets and hyperbolicity of our maps
by using properties of our maps and Kobayashi metrics.

%%%%%%%%%%%%%%%%%%%%%%%%%%%%%%%%%%%%%%%%%%%%%%%%%%%%%%%%%%%%%%%%%%%%%%%%%%%%%%%%%%%%%%%%%%%
\section{$S_{k+2}$-equivariant maps}

Crass~\cite{c} selected the symmetric group $S_{k+2}$ as a finite group
acting on $\mathbf{P}^{k}$ and found an $S_{k+2}$-\hspace{0pt}\textit{equivariant} map
which is holomorphic and \textit{critically finite} for each $k \ge 1$.
We denote by $C=C(f)$ the critical set of $f$ and say that $f$ is \textit{critically finite}
if each irreducible component of $C(f)$ is periodic or preperiodic.
More precisely, $S_{k+2}$-\textit{equivariant} map $g_{k+3}$
defined in section \ref{formula} preserves
each irreducible component of $C(g_{k+3})$, which is a projective hyperplane.
The complement of $C(g_{k+3})$ is Kobayashi hyperbolic.
Furthermore restrictions of $g_{k+3}$ to invariant projective subspaces
have the same properties as above.
See section \ref{properties} for details.

%%%%%%%%%%%%%%%%%%%%%%%%%%%%%%%%%%%%%%%%%%%%%%%%%%%%%%%%%%%%%%%%%%%%%%%%%%%%%%%%%%%%%%%%%%%
\subsection{$S_{k+2}$ acts on  $\mathbf{P}^{k}$}

An action of the $(k+2)$-th symmetric group $S_{k+2}$ on $\mathbf{P}^{k}$ is induced by
the permutation action of $S_{k+2}$ on $\mathbf{C}^{k+2}$
for each $k \ge 1$.
The transposition $(i,j)$ in $S_{k+2}$ corresponds with
the transposition $"u_i \leftrightarrow u_j"$  on $\mathbf{C}^{k+2}_u$,
which pointwise fixes the hyperplane $\{ u_i = u_j \}
= \{ u \in \mathbf{C}^{k+2}_u \ | \ u_i = u_j \}$.
Here $\mathbf{C}^{k+2} = \mathbf{C}^{k+2}_u = \{ u=(u_1,u_2,\cdot \cdot,u_{k+2})
\ | \ u_i \in \mathbf{C} \ \ \text{for} \ \ i=1,\cdot \cdot,k+2 \}$.

The action of $S_{k+2}$ preserves a hyperplane $H$ in $\mathbf{C}^{k+2}_u$,
which is identified with $\mathbf{C}^{k+1}_x$
by projection $A:\mathbf{C}^{k+2}_u  \to \mathbf{C}^{k+1}_x$,
\[
H= \left\{ \sum_{i=1}^{k+2} u_i =0 \right\}
\stackrel{\mathrm{A}}{\simeq} \mathbf{C}^{k+1}_{x}
\text{ and } \ A=
\left(
\begin{array}{ccccc}
1 & 0 & \ldots & 0 & -1 \\
0 & 1 & \ldots & 0 & -1 \\
\vdots & \vdots & \ddots & \vdots & \vdots \\
0 & 0 & \ldots & 1 & -1
\end{array}
\right).
\]
Here $\mathbf{C}^{k+1} = \mathbf{C}^{k+1}_x =\{ x=(x_1,x_2,\cdot \cdot,x_{k+1})
\ | \ x_i \in \mathbf{C} \ \ \text{for} \ \ i=1,\cdot \cdot,k+1 \}$.

Thus the permutation action of $S_{k+2}$ on $\mathbf{C}^{k+2}_u$
induces an action of $"S_{k+2}"$ on $\mathbf{C}^{k+1}_x$.
Here $"S_{k+2}"$ is generated by the permutation action $S_{k+1}$ on $\mathbf{C}^{k+1}_x$
and a $(k+1,k+1)$-matrix $T$ which corresponds to the transposition $(1 ,\ k+2)$ in $S_{k+2}$,
\[
T=
\left(
\begin{array}{cccc}
-1 & 0 & \ldots & 0 \\
-1 & 1 & \ldots & 0 \\
\vdots & \vdots & \ddots & 0 \\
-1 & 0 & \ldots & 1
\end{array}
\right).
\]
Hence the hyperplane corresponding to $\{ u_i = u_j \}$
is $\{ x_i = x_j \}$ for $1 \leq i < j \leq k+1$.
The hyperplane corresponding to $\{ u_i = u_{k+2} \}$
is $\{ x_i =0 \}$ for $1 \leq i \leq k+1$.
Each element in $"S_{k+2}"$ which corresponds to some transposition in $S_{k+2}$
pointwise fixes one of these hyperplanes in $\mathbf{C}^{k+1}_x$.

The action of $"S_{k+2}"$ on $\mathbf{C}^{k+1}$
projects naturally to the action of $"S_{k+2}"$ on $\mathbf{P}^{k}$.
These hyperplanes on $\mathbf{C}^{k+1}$
projects naturally to projective hyperplanes on $\mathbf{P}^{k}$.
Here $\mathbf{P}^{k} =
\{ x=[x_1:x_2:\cdot \cdot:x_{k+1}] \ | \ (x_1,x_2,\cdot \cdot,x_{k+1})
\in \mathbf{C}^{k+1} \setminus \{ \mathbf{0} \} \}$.
Each element in the action of $"S_{k+2}"$ on $\mathbf{P}^{k}$
which corresponds to some transposition in $S_{k+2}$
pointwise fixes one of these projective hyperplanes.
We denote $"S_{k+2}"$ also by $S_{k+2}$
and call these projective hyperplanes transposition hyperplanes.

%%%%%%%%%%%%%%%%%%%%%%%%%%%%%%%%%%%%%%%%%%%%%%%%%%%%%%%%%%%%%%%%%%%%%%%%%%%%%
\subsection{Existence of our maps}
\label{formula}

One way to get $S_{k+2}$-\textit{equivariant} maps on $\mathbf{P}^{k}$
which are \textit{critically finite} is to make $S_{k+2}$-\textit{equivariant} maps
whose critical sets coincide with the union of the transposition hyperplanes.

\begin{thm}[\cite{c}]\label{c}
For each $k \ge 1$,
$g_{k+3}$ defined below is the unique $S_{k+2}$-\textit{equivariant} holomorphic map
of degree $k+3$ which is doubly critical on each transposition hyperplane.
\[
g=g_{k+3} =[g_{k+3,1}:g_{k+3,2}:\cdot \cdot:g_{k+3,k+1}]:\mathbf{P}^{k} \to \mathbf{P}^{k},
\]
\[
\text{where } \ g_{k+3,l} (x) =x_l^3 \sum_{s=0}^{k} (-1)^s \frac{s+1}{s+3} x_l^s A_{k-s},
\ \ A_0 =1,
\]
\[
\text{ and } A_{k-s} \text{ is the elementary symmetric function of degree k-s in }
\mathbf{C}^{k+1}.
\]
\end{thm}

Then the critical set of $g$ coincides with the union of the transposition hyperplanes.
Since $g$ is $S_{k+2}$-\hspace{0pt}\textit{equivariant}
and each transposition hyperplane is pointwise fixed by some element in $S_{k+2}$,
$g$ preserves each transposition hyperplane.
In particular $g$ is \textit{critically finite}.
Although Crass~\cite{c} used this explicit formula
to prove Theorem {\rmfamily \ref{c}},
we shall only use properties of the $S_{k+2}$-\hspace{0pt}\textit{equivariant} maps
described below.

%%%%%%%%%%%%%%%%%%%%%%%%%%%%%%%%%%%%%%%%%%%%%%%%%%%%%%%%%%%%%%%%%%%%%%%%%%%%%%%%%%%%
\subsection{Properties of our maps}
\label{properties}

Let us look at properties of the $S_{k+2}$-\textit{equivariant} map $g$ on $\mathbf{P}^{k}$
for a fixed $k$, which is proved in \cite{c} and shall be used to prove our results.
Let $L^{k-1}$ denote one of the transposition hyperplanes,
which is isomorphic to $\mathbf{P}^{k-1}$.
Let $L^m$ denote one of the intersections of $(k-m)$ or more distinct transposition hyperplanes
which is isomorphic to $\mathbf{P}^{m}$ for $m=0,1, \cdot \cdot ,k-1$.

First, let us look at properties of $g$ itself.
The critical set of $g$ consists of the union of the transposition hyperplanes.
By $S_{k+2}$-\hspace{0pt}\textit{equivariance}, $g$ preserves each transposition hyperplane.
Furthermore the complement of the critical set of $g$ is Kobayashi hyperbolic.

Next, let us look at properties of $g$ restricted to $L^m$ for $m=1,2, \cdot \cdot ,k-1$.
Let us fix any $m$. Since $g$ preserves each $L^m$,
we can also consider the dynamics of $g$ restricted to any $L^m$.
Each restricted map has the same properties as above.
Let us fix any $L^m$ and denote by $g|_{L^m}$ the restricted map of $g$ to the $L^m$.
The critical set of $g|_{L^m}$ consists of the union of 
intersections of the $L^m$ and another $L^{k-1}$ which does not include the $L^m$.
We denote it by $L^{m-1}$,
which is an irreducible component of the critical set of $g|_{L^m}$.
By $S_{k+2}$-\hspace{0pt}\textit{equivariance},
$g|_{L^m}$ preserves each irreducible component of the critical set of $g|_{L^m}$.
Furthermore the complement of the critical set of $g|_{L^m}$ in $L^m$ is Kobayashi hyperbolic.

Finally, let us look at a property of superattracting fixed points of $g$.
The set of superattracting points,
where the derivative of $g$ vanishes for all directions,
coincides with the set of $L^0$'s.

\begin{rem}\label{} 
For every $k \geq 1$ and every $m$, $1 \leq m \leq k$,
a restricted map of $g_{k+3}$ to any $L^m$ is not conjugate to $g_{m+3}$.
\end{rem}

%%%%%%%%%%%%%%%%%%%%%%%%%%%%%%%%%%%%%%%%%%%%%%%%%%%%%%%%%%%%%%%%%%%%%%%%%%%%%%%%%%%%%%%%
\subsection{Examples for $k=1$ and $2$}

Let us see transposition hyperplanes of the $S_3$-\textit{equivariant} function $g_4$
and the $S_4$-\textit{equivariant} map $g_5$ to make clear what $L^m$ is.
In \cite{c} one can find explicit formulas and figures of dynamics of
$S_{k+2}$-equivariant maps in low-dimensions .

\subsubsection{$S_3$-equivariant function $g_4$ in $\mathbf{P}^{1}$}

\[
g_3([x_1:x_2])=[x_1^3 (-x_1+2x_2):x_2^3 (2x_1-x_2)]
:\mathbf{P}^{1} \to \mathbf{P}^{1},
\]
\[
C(g_3)= \{x_1=0 \} \cup \{x_2=0 \} \cup \{x_1=x_2 \} = \{0,1, \infty \}
\text{ in } \ \mathbf{P}^{1} .
\]
In this case "transposition hyperplanes" are points in $\mathbf{P}^{1}$ and
$L^0$ denotes one of three superattracting fixed points of $g_3$.

\subsubsection{$S_4$-equivariant map $g_5$ in $\mathbf{P}^{2}$}

\[
C(g_5)= \{x_1=0 \} \cup \{x_2=0 \} \cup \{x_3=0 \} \cup 
\]
\[
\{x_1=x_2 \} \cup \{x_2=x_3 \} \cup \{x_3=x_1 \}
\ \text{in} \ \mathbf{P}^{2}.
\]
In this case $L^1$ denotes one of six transposition hyperplanes in $\mathbf{P}^{2}$,
which is an irreducible component of $C(g_5)$.
For example, let us fix a transposition hyperplane $\{x_1=0 \}$.
Since $g_5$ preserves each transposition hyperplane, we can also consider
the dynamics of $g_5$ restricted to $\{x_1=0 \}$.
We denote by $g_5 |_{ \{x_1=0 \} }$ the restricted map of $g_5$ to $\{x_1=0 \}$.
The critical set of $g_5 |_{ \{x_1=0 \} }$ in $\{x_1=0 \} \simeq  \mathbf{P}^{1}$ is
\[
C(g_5 |_{ \{x_1=0 \} }) = \{ [0:1:0],[0:0:1],[0:1:1] \}.
\]
When we use $L^0$ after we fix $\{x_1=0 \}$, $L^0$ denotes one of intersections of
$\{x_1=0 \}$ and another transposition hyperplane,
which is a superattracting fixed point of $g_5 |_{ \{x_1=0 \} }$ in $\mathbf{P}^{1}$.
The set of superattracting fixed points of $g_5$ in $\mathbf{P}^{2}$ is
\[
\{ [1:0:0],[0:1:0],[0:0:1],[1:1:1],[1:1:0],[1:0:1],[0:1:1] \}.
\]
In general $L^0$ denotes one of intersections of two or more transposition hyperplanes,
which is a superattracting fixed point of $g_5$ in $\mathbf{P}^{2}$.

%%%%%%%%%%%%%%%%%%%%%%%%%%%%%%%%%%%%%%%%%%%%%%%%%%%%%%%%%%%%%%%%%%%%%%%%%%%%%%%%%%%%%%%%%%
\section{The Fatou sets of the $S_{k+2}$-equivariant maps} 

%%%%%%%%%%%%%%%%%%%%%%%%%%%%%%%%%%%%%%%%%%%%%%%%%%%%%%%%%%%%%%%%%%%%%%%%%%%%%%%%%%%%%%%%%%%
\subsection{Definitions and preliminaries}

Let us recall theorems about \textit{critically finite} holomorphic maps.
Let $f$ be a holomorphic map from $\mathbf{P}^{k}$ to $\mathbf{P}^{k}$.
The Fatou set of $f$ is defined to be the maximal open subset
where the iterates $\{ f^n \}_{n \geq 0}$ is a normal family.
The Julia set of $f$ is defined to be the complement of the Fatou set of $f$.
Each connected component of the Fatou set is called a Fatou component.
Let $U$ be a Fatou component of $f$.
A holomorphic map $h$ is said to be a limit map on $U$
if there is a subsequence $\{ f^{n_s} |_U \}_{s \geq 0}$
which locally converges to $h$ on $U$.
We say that a point $q$ is a Fatou limit point
if there is a limit map $h$ on a Fatou component $U$ such that $q \in h(U)$.
The set of all Fatou limit points is called the Fatou limit set.
We define the $\omega$-limit set $E(f)$ of the critical points by
\[
E(f)= \bigcap_{j=1}^{\infty} \overline{ \bigcup^{\infty}_{n=j} f^n(C) }.
\]

\begin{thm}(\cite[Proposition 5.1]{ueda})\label{u}
If $f$ is a critically finite holomorphic map from $\mathbf{P}^{k}$ to $\mathbf{P}^{k}$,
then the Fatou limit set is contained in the $\omega$-limit set $E(f)$.
\end{thm}

Let us recall the notion of Kobayashi metrics.
Let $M$ be a complex manifold 
and $K_M (x,v)$ the Kobayashi quasimetric on $M$,
\[
\inf \left\{ |a| \Bigl| \varphi : \mathbf{D} \to M : \text{holomorphic}
,\varphi (0)=x
,D{\varphi} \left(a \left( \frac{ \partial }{ \partial z} \right)_0 \right)=v
,a \in \mathbf{C} \right\}
\]
for $x \in M,\ v \in T_x M,\ z \in \mathbf{D}$,
where $\mathbf{D}$ is the unit disk in $\mathbf{C}$.
We say that $M$ is Kobayashi hyperbolic if $K_M$ becomes a  metric.
Theorem {\rmfamily \ref{u1}} is a corollary of
Theorem {\rmfamily \ref{m}} and Theorem {\rmfamily \ref{fs}}
for $k=1$ and $2$.

\begin{thm}\label{m}(a basic result whose former statement can be found in \cite[Corollary 14.5]{m})
If $f$ is a critically finite holomorphic function from $\mathbf{P}^{1}$ to $\mathbf{P}^{1}$,
then the only Fatou components of $f$ are attractive components of superattracting points.
Moreover if the Fatou set is not empty, then the Fatou set has full measure in $\mathbf{P}^{1}$.
\end{thm}

\begin{thm}\label{fs}(\cite[theorem 7.7]{fs1})
If $f$ is a critically finite holomorphic map from $\mathbf{P}^{2}$ to $\mathbf{P}^{2}$
and the complement of $C(f)$ is Kobayashi hyperbolic,
then the only Fatou components of $f$ are attractive components of superattracting points.
\end{thm}

%%%%%%%%%%%%%%%%%%%%%%%%%%%%%%%%%%%%%%%%%%%%%%%%%%%%%%%%%%%%%%%%%%%%%%%%%%%%%%%%%%%
\subsection{Our first result}

%%%%%%%%%%%%%%%%%%%%%%%%%%%%%%%%%%%%%%%%%%%%%%%%%%%%%%%%%%%%%%%%%%%%%%%%%%%%%%%%%%%%%%%%%%%%%
Let us fix any $k$ and $g=g_{k+3}$.
For every $m$, $2 \leq m \leq k$,
we can apply an argument in \cite{fs1} to a restricted map of $g$ to any $L^{m}$
because every $L^{m-1}$ is smooth
and because every $L^m \backslash C(g|_{L^m})$ is Kobayashi hyperbolic.
We shall use this argument in Lemma {\rmfamily \ref{lem}},
which is used to prove Proposition {\rmfamily \ref{pro}}.

%%%%%%%%%%%%%%%%%%%%%%%%%%%%%%%%%%%%%%%%%%%%%%%%%%%%%%%%%%%%%%%%%%%%%%%%%%%%%%%%%%%%%%%%%%%%%
\begin{pro}\label{pro}
For any Fatou component $U$ which is disjoint from $C(g)$,
there exists an integer $n$ such that $g^{n}(U)$ intersects with $C(g)$.
\end{pro}

\begin{proof}[Proof:]
We suppose that $g^n (U)$ is disjoint from $C(g)$ for any $n$
and derive a contradiction by using Lemma {\rmfamily \ref{lem}}
and Remark {\rmfamily \ref{l.rem}} below.
Take any point $x_{0} \in U$. Since $E(g)$ coincides with $C(g)$,
$g^n(x_{0})$ accumulates to $C(g)$ as $n$ tends to $\infty$ from Theorem {\rmfamily \ref{u}}.
Since $C(g)$ is the union of the transposition hyperplanes,
there exists a smallest integer $m_1$ such that $g^n(x_{0})$ accumulates to some $L^{m_1}$.
Let $h_1$ be a limit map on $U$ such that $h_1 (x_{0})$ belongs to the $L^{m_1}$.
From Lemma {\rmfamily \ref{lem}} below,
the intersection of $h_1 (U)$ and the $L^{m_1}$ is an open set in the $L^{m_1}$
and is contained in the Fatou set of $g|_{L^{m_1}}$.

We next consider the dynamics of $g|_{L^{m_1}}$.
If there exists an integer $n_2$
such that $g^{n_{2}}(h_1 (U) \cap L^{m_1})$ intersects with $C(g|_{L^{m_1}})$,
then $g^{n_{2}}(h_1 (U) \cap L^{m_1})$ intersects with some $L^{m_1 -1}$.
In this case we can consider the dynamics of $g|_{L^{m_1 -1}}$.
On the other hand, if there does not exist such $n_2$,
then there exists an integer $m_2$ and a limit map $h_2$ on $h_1 (U) \cap L^{m_1}$
such that the intersection of $h_2 (h_1 (U) \cap L^{m_1})$ and some $L^{m_2}$ 
is an open set in the $L^{m_2}$ from Remark {\rmfamily \ref{l.rem}} below.
Thus it is contained in the Fatou set of $g|_{L^{m_2}}$.
Here $m_2$ is smaller than $m_1$.
In this case we can consider the dynamics of $g|_{L^{m_2}}$.

We continue the same argument above.
These reductions finally come to some $L^1$ and we use Theorem {\rmfamily \ref{m}}.
One can find a similar reduction argument
in the proof of Theorem {\rmfamily \ref{u1}}.
Consequently $g^n(x_{0})$ accumulates to some superattracting point $L^0$.
So there exists an integer $s$ such that $g^s$ sends $U$ to
the attractive Fatou component which contains the superattracting point $L^0$.
Thus $g^{s}(U)$ intersects with $C(g)$, which is a contradiction.
\end{proof}

%%%%%%%%%%%%%%%%%%%%%%%%%%%%%%%%%%%%%%%%%%%%%%%%%%%%%%%%%%%%%%%%%%%%%%%%%%%%%%%%%%%%%%%%%%%%%
\begin{rem}\label{p.rem}
Even if a Fatou component $U$ intersects with some $L^m$
and is disjoint from any $L^{m-1}$,
then the similar thing as above holds for the dynamics in the $L^m$.
In this case
$U \cap L^m$ is contained in the Fatou set of $g|_{L^m}$
and there exists an integer $n$ such that
$g^{n}(U \cap L^m)$ intersects with $C(g|_{L^m})$.
\end{rem}

%%%%%%%%%%%%%%%%%%%%%%%%%%%%%%%%%%%%%%%%%%%%%%%%%%%%%%%%%%%%%%%%%%%%%%%%%%%%%%%%%%%%%%%%
\begin{lem}\label{lem}
For any Fatou component $U$ which is disjoint from $C(g)$ and any point $x_{0} \in U$,
let $h$ be a limit map on $U$ such that 
$h(x_{0})$ belongs to some $L^m$ and does not belong to any $L^{m-1}$.
If $g^n(U)$ is disjoint from $C(g)$ for every $n \ge 1$,
then the intersection of $h(U)$ and the $L^m$ is an open set in the $L^m$.
\end{lem}

\begin{proof}[Proof:]
Let $B$ be the complement of $C(g)$.
Since $B$ is Kobayashi hyperbolic and $B$ includes $g^{-1}(B)$,
$g^{-1}(B)$ is Kobayashi hyperbolic, too.
So we can use Kobayashi metrics $K_B$ and $K_{g^{-1}(B)}$.
Since $B$ includes $g^{-1}(B)$,
\[
K_B (x,v) \le K_{g^{-1}(B)} (x,v)
\ \text{ for all $x \in g^{-1}(B)$,\ $v \in T_x \mathbf{P}^k$ }.
\]
In addition, since $g$ is an unbranched covering from $g^{-1}(B)$ to $B$,
\[
K_{g^{-1}(B)} (x,v) = K_B (g(x),Dg(v))
\ \text{ for all $x \in g^{-1}(B)$,\ $v \in T_x \mathbf{P}^k$ }.
\]
From these two inequalities we have the following inequality
\[
K_B (x,v) \le K_B (g(x),Dg(v)) \ \text{ for all $x \in g^{-1}(B)$,\ $v \in T_x \mathbf{P}^k$ }.
\]
Since the same argument holds for any $g^n$
from $g^{-n}(B)$ to $B$,
\[
K_B (x,v) \le K_B (g^n (x),Dg^n (v))
\ \text{ for all $x \in g^{-n}(B)$,\  $v \in T_x \mathbf{P}^k$ }.
\]
Since $g^n$ is an unbranched covering from $U$ to $g^n (U)$ and
$B$ includes $g^n (U)$ for every $n$,
a sequence $\{ K_B (g^n (x),Dg^n (v)) \}_{n \geq 0}$
is bounded for all $x \in U$, $v \in T_x \mathbf{P}^k$.
Hence we have the following inequality for any unit vectors $v_n$ in $T_{x_{0}} U$
with respect to the Fubini-Study metric in $\mathbf{P}^k$,
\begin{equation}
0< \inf_{|v|=1} K_B (x_{0},v) \le K_B (x_{0},v_n) \le K_B (g^n (x_{0}),Dg^n (x_{0})v_n) < \infty .
\end{equation}
That is, the sequence $\{ K_B (g^n (x_{0}),Dg^n (x_{0})v_n) \}_{n \geq 0}$
is bounded away from $0$ and $\infty$ uniformly.

We shall choose $v_n$ so that $Dg^n (x_{0})v_n$ keeps parallel to the $L^m$
and claim that $Dh(x_{0})v \neq \mathbf{0}$
for any accumulation vector $v$ of $v_n$.
Let $h= \lim_{n \to \infty} g^n$ for simplicity.
Let $V$ be a neighborhood of $h(x_0)$ and $\psi$ a local coordinate on $V$
so that $\psi (h(x_{0}))= \mathbf{0}$ and $\psi (L^m \cap V) \subset 
\{ y=(y_1,y_2, \cdot \cdot ,y_k) \ | \ y_1 = \cdot \cdot = y_{k-m} =0 \}$.
In this chart there exists a constant $r > 0$ such that
a polydisk $P(\mathbf{0},2r)$ does not intersect with
any images of transposition hyperplanes which do not include the $L^m$.
Since $\psi (g^n (x_{0}))$ converges to $\mathbf{0}$ as $n$ tends to $\infty$,
we may assume that $\psi (g^n (x_{0}))$ belongs to $P(\mathbf{0},r)$ for large $n$.
Let $\{ v_n \}_{n \geq 0}$ be unit vectors in $T_{x_{0}} \mathbf{P}^k$
and $\{ w_n \}_{n \geq 0}$ vectors in $T_{\psi (g^n(x_0))} \mathbf{C}^k$
so that $w_n$ keep parallel to $\psi (L^m)$ with a same direction
and 
\[
Dg^n (x_{0})v_n =| Dg^n (x_{0})v_n | \ D \psi^{-1} (w_n).
\]
So we may assume that the length of $w_n$ is almost unit for large $n$.
We define holomorphic maps $\varphi_{n}$ from $\mathbf{D}$ to $P(\mathbf{0},2r)$ as
\[
\varphi_{n} (z) = \psi (g^n (x_{0})) +rzw_n \ \text{ for } z \in \mathbf{D}
\]
and consider holomorphic maps $\psi^{-1} \circ \varphi_{n}$ from $\mathbf{D}$ to $B$
for large $n$.
Then
\[
(\psi^{-1} \circ \varphi_{n})(0) = g^n (x_{0}),
\]
\[
D(\psi^{-1} \circ \varphi_{n})
\left( \frac{| Dg^n (x_{0})v_n |}{r} \left( \frac{\partial}{\partial z} \right)_0 \right)
=Dg^n (x_{0})v_n .
\]
Suppose $Dh(x_{0})v = \mathbf{0}$, then $Dg^n (x_{0})v$ converges to $\mathbf{0}$
as $n$ tends to $\infty$ and so does $Dg^n (x_{0})v_n$.
By the definition of Kobayashi metric we have that
\[
K_B (g^n (x_{0}),Dg^n (x_{0})v_n) \le \frac{| Dg^n (x_{0})v_n |}{r} \to 0
\ \text{ as } \ n \to \infty .
\]
Since this contradicts (1), we have $Dh(x_{0})v \neq \mathbf{0}$.
This holds for all directions which are parallel to $\psi (L^m)$.
Consequently the intersection of $h(U)$ and the $L^m$ is an open set in $L^m$.
\end{proof}

%%%%%%%%%%%%%%%%%%%%%%%%%%%%%%%%%%%%%%%%%%%%%%%%%%%%%%%%%%%%%%%%%%%%%%%%%%%%%%%%%%%%%%%%%%%%%
\begin{rem}\label{l.rem}
The similar thing as above holds for the dynamics of any restricted map.
Thus even if a Fatou component $g^n(U)$ intersects with $C(g)$ for some $n$,
the same result as above holds.
Because one can consider the dynamics in the $L^m$
when $g^n(U)$ intersects with some $L^m$.
\end{rem}

%%%%%%%%%%%%%%%%%%%%%%%%%%%%%%%%%%%%%%%%%%%%%%%%%%%%%%%%%%%%%%%%%%%%%%%%%%%%%%%%%%%%%%%%%%%%%
\begin{thm}\label{u1}
For each $k \ge 1$,
the Fatou set of the $S_{k+2}$-equivariant map $g$ consists of
attractive basins of superattracting fixed points
which are intersections of $k$ or more distinct transposition hyperplanes.
\end{thm}

\begin{proof}[Proof:]
This theorem follows from Proposition {\rmfamily \ref{pro}}
and Remark {\rmfamily \ref{p.rem}} immediately.
Let us describe details.
Take any Fatou component $U$.
From Proposition {\rmfamily \ref{pro}}
there exists an integer $n_k$ such that $g^{n_{k}}(U)$ intersects with $C(g)$.
Since $C(g)$ is the union of the transposition hyperplanes,
$g^{n_{k}}(U)$ intersects with some $L^{k-1}$.
By doing the same thing as above for the dynamics of $g$ restricted to the $L^{k-1}$,
there exists an integer $n_{k-1}$ such that
$g^{n_k + n_{k-1}}(U)$ intersects with some $L^{k-2}$
from Remark {\rmfamily \ref{p.rem}}.
We again do the same thing as above for the dynamics of $g$ restricted to the $L^{k-2}$.

These reductions finally come to some $L^1$.
That is, there exists integers $n_{k-2}, \cdot \cdot ,n_{2}$ such that
$g^{n_k + n_{k-1} + \cdot \cdot + n_2} (U)$ intersects with some $L^{1}$.
From Theorem {\rmfamily \ref{m}} there exists an integer $n_1$ such that
$g^{n_1}(g^{n_k + n_{k-1} + \cdot \cdot + n_2} (U))$ contains some $L^0$.
Hence $g^{n_k + n_{k-1} + \cdot \cdot + n_1}$ sends $U$ to the attractive Fatou component
which contains the superattracting fixed point $L^0$ in $\mathbf{P}^{k}$.
\end{proof}

%%%%%%%%%%%%%%%%%%%%%%%%%%%%%%%%%%%%%%%%%%%%%%%%%%%%%%%%%%%%%%%%%%%%%%%%%%%%%%%%%%%%%%%%%%%
\section{Axiom A and the $S_{k+2}$-equivariant maps}

%%%%%%%%%%%%%%%%%%%%%%%%%%%%%%%%%%%%%%%%%%%%%%%%%%%%%%%%%%%%%%%%%%%%%%%%%%%%%%%%%%%%%%%%%%%
\subsection{Definitions and preliminaries}

Let us define hyperbolicity of non-invertible maps and the notion of Axiom A.
See \cite{mj} for details.
Let $f$ be a holomorphic map from $\mathbf{P}^k$ to $\mathbf{P}^k$
and $K$ a compact subset such that $f( K ) = K$.
Let $\widehat{K}$ be the set of histories in $K$ and
$\widehat{f}$ the induced homeomorphism on $\widehat{K}$.
We say that $f$ is hyperbolic on $K$ if there exists a continuous decomposition
$T_{\widehat{K}} =E^u + E^s$ of the tangent bundle such that
$D \widehat{f} (E_{\widehat{x}}^{u/s}) \subset E_{\widehat{f}(\widehat{x})}^{u/s}$  
and if there exists constants $c>0$ and $\lambda >1$ such that for every $n \ge 1$,
\[
| D \widehat{f}^n (v) | \ge c \lambda^n |v| \ \text{ for all } v \in E^u \ \text{and}
\]
\[
| D \widehat{f}^n (v) | \le c^{-1} \lambda^{-n} |v| \ \text{ for all } v \in E^s.
\]
Here $| \cdot |$ denotes the Fubini-Study metric on $\mathbf{P}^k$.
If a decomposition and inequalities above hold for $f$ and $K$,
then it also holds for $\widehat{f}$ and $\widehat{K}$.
In particular we say that $f$ is expanding on $K$ if $f$ is hyperbolic on $K$
with unstable dimension $k$.
Let $\Omega$ be the non-wandering set of $f$, i.e.,
the set of points for any neighborhood $U$ of which there exists an integer $n$
such that $f^n (U)$ intersects with $U$.
By definition, $\Omega$ is compact and $f( \Omega )= \Omega$.
We say that $f$ satisfies Axiom A if $f$ is hyperbolic on $\Omega$ and
periodic points are dense in $\Omega$.

Let us introduce a theorem which deals with repelling part of dynamics.
Let $f$ be a holomorphic map from $\mathbf{P}^k$ to $\mathbf{P}^k$.
We define the k-th Julia set $J_k$ of $f$
to be the support of the measure with maximal entropy,
in which repelling periodic points are dense.
It is a fundamental fact that in dimension 1 the 1st Julia set $J_1$ coincides with the Julia set $J$.
Let $K$ be a compact subset such that $f( K ) = K$.
We say that $K$ is a repeller if $f$ is expanding on $K$.  

\begin{thm}\label{km1}(\cite{km})
Let $f$ be a holomorphic map on $\mathbf{P}^k$ of degree at least 2
such that the $\omega$-limit set $E(f)$ is pluripolar. 
Then any repeller for $f$ is contained in $J_k$. In particular,
\[
J_k = \overline{ \{ repelling \ periodic \ points \ of \ f \} }
\]
\end{thm}

If $f$ is critically finite, then $E(f)$ is pluripolar.
We need the theorem above to prove our second result.

%%%%%%%%%%%%%%%%%%%%%%%%%%%%%%%%%%%%%%%%%%%%%%%%%%%%%%%%%%%%%%%%%%%%%%%%%%%%%%%%%%%
\subsection{Our second result}

\begin{thm}\label{u2}
For each $k \ge 1$, the $S_{k+2}$-equivariant map $g$ satisfies Axiom A.
\end{thm}

\begin{proof}[Proof:]
We only need to consider 
the $S_{k+2}$-\hspace{0pt}\textit{equivariant} map $g$ for a fixed $k$,
because argument for any $k$ is similar as the following one.
Let us show the statement above for a fixed $k$ by induction.
A restricted map of $g$ to any ${L^1}$ satisfies Axiom A
by using the theorem of \textit{critically finite} functions (see \cite[Theorem 19.1]{m}).
We only need to show that a restricted map of $g$ to a fixed $L^2$ satisfies Axiom A.
Then a restricted map of $g$ to any $L^2$ satisfies Axiom A by symmetry.
Argument for a restricted map of $g$ to any $L^m$,
$3 \le m \le k$, is similar as for a restricted map of $g$ to the $L^2$.
Let us denote $g |_{L^2}$, $\Omega (g |_{L^2})$, and $L^2$
by $g$, $\Omega$, and $\mathbf{P}^{2}$ for simplicity.

We want to show that $g |_{L^2}$ is hyperbolic on $\Omega (g |_{L^2})$
by using Kobayashi metrics.
If $g$ is hyperbolic on $\Omega$,
then $\Omega$ has a decomposition to $S_i$,
\[
\Omega = S_0 \cup S_1 \cup S_2,
\]
where i=0,1,2 indicate the unstable dimensions.
Since $C(g)$ attracts all nearby points,
$S_0$ includes all the $L^0$'s and $S_1$ includes all the Julia sets of $g|_{L^1}$.
We denote by ${J(g|_{L^1})}$ the Julia set of $g|_{L^1}$.
Then $g$ is contracting in all directions at ${L^0}$ and
is contracting in the normal direction and expanding in an $L^1$-direction
on ${J(g|_{L^1})}$.
Let us consider a compact, completely invariant subset in $\mathbf{P}^{2} \setminus C$,
\[
S= \{ x \in \mathbf{P}^{2} \ | \ \text{dist} (g^n (x),C)
\nrightarrow 0 \ \text{ as } \ n \to \infty \}.
\]
By definition, we have $J_2 \subset S_2 \subset S$.
If $g$ is expanding on $S$,
then it follow that $S_0 = \cup L^0,\ S_1 = \cup J(g|_{L^1})$.
Moreover $J_2 = S_2 = S$ holds from Theorem {\rmfamily \ref{km1}}
(see Remark {\rmfamily \ref{km3}} below).
Since periodic points are dense in $J(g|_{L^1})$ and $J_2$,
expansion of $g$ on $S$ implies Axiom A of $g$.

Let us show that $g$ is expanding on $S$.
Because $f$ is attracting on $C$ and preserves $C$,
there exists a neighborhood $V$ of $C$ such that
$V$ is relatively compact in $g^{-1} (V)$
and the complement of $V$ is connected.  %% \mathbf{P}^{2} \setminus 
We assume one of $L^1$'s to be the line at infinity of $\mathbf{P}^{2}$.
By letting $B$ be $\mathbf{P}^{2} \setminus V$ and
$U$ one of connected components of $g^{-1}( \mathbf{P}^{2} \setminus V)$,
we have the following inclusion relations,
\[
U \subset g^{-1} (B) \Subset B \subset
\mathbf{C}^{2} = \mathbf{P}^{2} \setminus L^1.
\]
Because $B$ and $U$ are in a local chart, there exists a constant $\rho <1$ such that
\[
K_B (x,v) \le \rho K_U (x,v) \ \text{ for all } x \in U,\ v \in T_x \mathbf{C}^{2}.
\]
In addition, since the map $g$ from $U$ to $B$ is an unbranched covering,
\[
K_U (x,v)=K_B (g(x),Dg(v)) \ \text{ for all } x \in U,\ v \in T_x \mathbf{C}^{2}.
\]
From these two inequalities we have the following inequality
\[
K_B (x,v) \le \rho K_B (g(x),Dg(v))
\ \text{ for all } x \in g^{-1} (B),\ v \in T_x \mathbf{C}^{2}.
\]
Since $g$ preserves $S$,
which is contained in $g^{-n} (B)$ for every $n \ge 1$, 
\[
K_B (x,v) \le \rho^n K_B (g^n (x),Dg^n (v))
\ \text{ for all } x \in S,\ v \in T_x \mathbf{C}^{2}.
\]
Consequently we have the following inequality for $\lambda = {\rho}^{-1} > 1$,
\[
K_B (g^n (x),Dg^n (v)) \ge \lambda^n K_B (x,v)
\ \text{ for all } x \in S,\ v \in T_x \mathbf{C}^{2}.
\]
Since $K_B (x,v)$ is upper semicontinuous and $|v|$ is continuous,
$K_B (x,v)$ and $|v|$ may be different only by a constant factor.
There exists $c>0$ such that
\[
|Dg^n (x)v| \ge c \lambda^n |v| \ \text{ for all } x \in S,\ v \in T_x \mathbf{C}^{2}.
\]
Thus $g$ is expanding on $S$ and satisfies Axiom A.
\end{proof}

%%%%%%%%%%%%%%%%%%%%%%%%%%%%%%%%%%%%%%%%%%%%%%%%%%%%%%%%%%%%%%%%%%%%%%%%%%%%%%%%%%%%%%%%%%%%%
\begin{rem}\label{km3}
Unlike the case when $k=1$, it does not seem obvious that $S$ being a repeller implies 
$J_k = S$ when $k \geq 2$. 
\end{rem}

%%%%%%%%%%%%%%%%%%%%%%%%%%%%%%%%%%%%%%%%%%%%%%%%%%%%%%%%%%%%%%%%%%%%%%%%%%%%%%%%%%%%%%%%%%%%%
\begin{rem}\label{}
From \cite[Theorem 4.11]{bowen} and \cite{qian},
it follows that the Fatou set of the $S_{k+2}$-equivariant map $g$ has full measure in $\mathbf{P}^{k}$
for each $k \ge 1$.
\end{rem}

%%%%%%%%%%%%%%%%%%%%%%%%%%%%%%%%%%%%%%%%%%%%%%%%%%%%%%%%%%%%%%%%%%%%%%%%%%%%%%%%%%%%%%%%%%%%%
\hspace{-7mm}
\textbf{Acknowledgments.}
\vspace{0mm}
\hspace{0mm}
I would like to thank Professor S. Ushiki and Doctor K. Maegawa for their useful advice.
Particularly in order to obtain our second result, Maegawa's suggestion to use Theorem {\rmfamily \ref{km1}} was helpful. 

%%%%%%%%%%%%%%%%%%%%%%%%%%%%%%%%%%%%%%%%%%%%%%%%%%%%%%%%%%%%%%%%%%%%%%%%%%%%%%%%%%%%%%%%%%%%%

%\begin{quote}
\hspace{-7mm} Graduate School of Human and Environmental Studies
\\ Kyoto University 
\\ Yoshida-Nihonmatu-cho, Sakyou-ku
\\ Kyoto 606-8501
\\ Japan
\\ \textit{E-mail address}: \texttt{ueno@math.h.kyoto-u.ac.jp}
%\end{quote}

\end{document}